\documentclass[onecolumn,12pt]{elsarticle}

\usepackage{amssymb}
\usepackage{amsfonts,amsmath,latexsym,amssymb}

\journal{arXiv.org}

\newtheorem{theorem}{Theorem}[section]
\newtheorem{definition}[theorem]{Definition}
\newtheorem{lemma}[theorem]{Lemma}
\newtheorem{remark}[theorem]{Remark}
\newtheorem{example}[theorem]{Example}
\newtheorem{corollary}[theorem]{Corollary}
\newtheorem{application}[theorem]{Application}

\def\Dj{\hbox{D\kern-.70em\raise.35ex\hbox{-} \raise-.10ex\hbox{}}}
\def\dj{\hbox{d\kern-.56em\raise.10ex\hbox{-} \raise-.10ex\hbox{\kern-.20em}}}

\begin{document}

\begin{frontmatter}

\title{\Large \bf Some considerations in relation to the \\
matrix equation \textit{\textbf{A}}\textit{\textbf{X}}\textit{\textbf{B}}$\,$=$\,$\textit{\textbf{C}}}

\medskip

\author{{\rm Biljana Radi\v ci\' c\mbox{\small ${\,}^{a}$}, Branko Male\v sevi\' c\mbox{\small ${\,}^{b}$}}}

\address[grf]{$\,${\rm Faculty of Civil Engineering, University of Belgrade, Serbia}}
\address[etf]{$\,${\rm Faculty of Electrical Engineering, University of Belgrade, Serbia}}

\begin{abstract}
{\small In this paper we represent a new form of condition for the consistency of the
matrix equation $AXB=C$. If the matrix equation $AXB=C$ is consistent, we determine a
form of general solution which contains both reproductive and non-reproductive
solutions. Also, we consider applications of the~concept~of~reproductivity for
obtaining general solutions of some matrix systems which are in relation to the matrix
equation $AXB=C$.}
\end{abstract}


\end{frontmatter}

\renewcommand{\thefootnote}{}

\footnotetext{$\!\!\!\!\!$ {\sl Email addresses}$\,:$

\smallskip
\mbox{\small $\!\!{\,}^{a\,}$}Biljana Radi\v ci\' c $<${\sl radicic.biljana@yahoo.com}$>$ (part time job) $\;$and$\;$
\mbox{\small $\!{\,}^{b\,}$}Branko Male\v sevi\' c $<${\sl malesevic@etf.rs}$>$}

\section{The reproductive equations}

\smallskip
The general concept of the reproductive equations was introduced by S.B. Pre\v si\' c
\cite{Presic68} in 1968. In this part of the paper we give the definition of
reproductive equations and the most important statements related to the reproductive
equations. Using the concept of reproductivity, in the next section, we obtain the
general solutions of some matrix systems which are in relation to the matrix
equation~$AXB=C$.

\smallskip
Let $S$ be a given non-empty set and $J$ be a given unary relation of $S$. Then
an~equation $J(x)$ is
{\em consistent} if there is at least one element $x_{0} \in S$, so-called {\em the
solution}, such~that $J(x_{0})$ is true. A formula $x=\phi(t)$,
where $\phi: S \longrightarrow S$ is a given function, represents
\mbox{\em the general solution} of the equation~$J(x)$ if~and~only~if
\vspace*{-1 mm}
$$
(\forall \, t) J(\phi(t)) \wedge (\forall \, x)( J(x) \Longrightarrow (\exists \, t) x
= \phi(t) ).
$$

\vspace*{-1 mm}

\noindent
In this part of the paper we give the definition of reproductive equations and the
fundamental statements related to the reproductive~equations.

\begin{definition}
\textit{The reproductive equations} are the equations of the following form:
\vspace*{-1 mm}
$$
x=\varphi(x),
$$

\vspace*{-1 mm}

\noindent
where $x$ is a unknown, $S$ is a given set and  $\varphi:S \longrightarrow S$ is a
given function which satisfies the following condition:
\vspace*{-1 mm}
\begin{equation}
\label{UR}
\varphi\circ\varphi=\varphi.
\end{equation}
\end{definition}

\break

\noindent
The condition (\ref{UR}) is called \textit{the condition of reproductivity}. The fundamental properties of the reproductive equations are
given by the following two statements S.B. Pre\v si\'c \cite{Presic68} (see also \cite{Bozic75},~\mbox{\cite{Presic72}-\cite{Presic00}}
and \cite{Tribute01}).

\vspace*{-0.5 mm}
\begin{theorem}
\label{T11}
For any consistent equation $J(x)$ there is an equation of the form $x=\varphi(x)$,
which is equivalent to $J(x)$ being in the same time reproductive as well.
\end{theorem}

\vspace*{-1.5 mm}

\begin{theorem}
\label{T12}
If a certain equation $J(x)$ is equivalent to the reproductive one $x=\varphi(x)$, the
general solution is given by the formula $x=\varphi(y)$, for any value $y\in S$.
\end{theorem}

\vspace*{-0.5 mm}
Let us remark that a formula $x=\phi(t)$, where $\phi \!:\!S \longrightarrow S$ is a
given function, represents \mbox{\em the reproductive general solution}
\cite{Bankovic11} of the equation~$J(x)$ if and only if
\vspace*{-0.5 mm}
$$
(\forall \, t) J(\phi(t)) \wedge (\forall \, t)( J(t) \Longrightarrow t = \phi(t) ).
$$

S.B. Pre\v si\' c  was the first one who considered implementations of reproductivity
on some matrix equations \cite{Presic68} (see also \cite{Haveric83}, \cite{Haveric84}
and \cite{Presic63}). The concept of reproductivity allows us to analyse various forms
of the solution. General applications of the concept of reproductivity were also
considered by J.D. Ke\v cki\' c in \cite{Keckic82}, \cite{Keckic83}, J.D. Ke\v cki\' c
and S.B. Pre\v si\' c in \cite{KeckicPresic97}, S. Rudeanu in
\cite{Rudeanu78}-\cite{Rudeanu01} and D. Bankovi\' c in
\cite{Bankovic79}-\cite{Bankovic11}.

\section{The matrix equation \textit{\textbf{A}}\!\,\textit{\textbf{X}}\!\,\textit{\textbf{B}}$\,$=$\,$\textit{\textbf{C}}}

\medskip
Let $m, n \in \mathbb{N}$ be natural numbers and $\mathbb{C}$ is the field of complex
numbers. The set of all matrices of order $ m \times n $ over $\mathbb{C}$ is
denoted by $\mathbb{C}^{m \times n}$. By $\mathbb{C}_{a}^{m \times n}$ we denote the
set of all $m \times n$ complex matrices of rank $a$. For $A \in
\mathbb{C}^{m \times n}$, the rank of $A$ is denoted by $\mbox{rank}(A)$. The unit
matrix of order $m$ is denoted by $I_m$ (if the dimension of unit matrix is known
from the context, we shall omit the index which indicates the dimension and use
the symbol $I$). Let $A=[a_{i,j}]\in \mathbb{C}^{m \times n}$. By $A_{i\rightarrow}$
and $A_{\downarrow j}$ we denote the $i$-th row of $A$ and the
$j$-th column of $A$, respectively. Therefore,
$$
A_{i \rightarrow }= ( a_{i,1}, a_{i,2},...,a_{i,n}),\,  i=1,...,m
$$
\vspace*{-2.5 mm}
and

$$
A_{ \downarrow j}=(a_{1,j},a_{2,j},...,a_{m,j})^{T}, \, j=1,...,n .
$$

\smallskip
\noindent
The matrix equation
\begin{equation}
\label{AXBC}
AXB=C
\end{equation}

\medskip
\noindent
was considered by many authors (\cite{Cvetkovic06}-\cite{Cvetkovic08},
\cite{Haveric83}, \cite{Haveric84}, \cite{Keckic85}-\cite{MalesevicRadicic12},
\cite{Presic63}, \cite{Presic00}, \cite{Tian10} and \cite{Tian12}). In the papers
\cite{ADajicJJKoliha}-\cite{Haveric84}, \cite{Keckic85} and \cite{Keckic97} the matrix
equation (\ref{AXBC}) was studied as a part of different matrix systems or as a
special case of corresponding matrix equations. Special case of the matrix equation
(\ref{AXBC}) is the following matrix equation:
\begin{equation}
\label{AXAA}
AXA=A.
\end{equation}
Any solution of this equation is called {\em $\{1\}$-inverse} of $A$ and
is~denoted~by~$A^{(1)}$. The set of all $\lbrace 1\rbrace$-inverses of $A$
is~denoted~by~$A\{1\}$.

\break

\noindent
For the matrix $A$, let regular matrices $Q \!\in\! \mathbb{C}^{m \times m}$ and
$P \!\in\! \mathbb{C}^{n \times n}$ be determined such that the following equality is
true:
\vspace*{-1.825 mm}
\begin{equation}
\label{EA}
QAP=E_{a}=
\left[
\begin{array}{c|c}
I_{a} & 0\\\hline
0 & 0
\end {array}
\right],
\end{equation}
where $a=\mbox{rank}(A)$. In \cite{Rohde64} C. Rohde showed that the general form of
$\{1\}$-inverse $A^{(1)}$ can be represented as:
\vspace*{-1.825 mm}
\begin{equation}
\label{AJEDAN}
A^{(1)} =
P
\left[
\begin{array}{c|c}
I_{a} & X_{1}\\\hline
X_{2} & X_{3}
\end {array}
\right]
Q \, ,
\end{equation}
where $X_{1}$, $X_{2}$ and $X_{3}$ are arbitrary matrices of suitable sizes (see also
\cite{Ben-IsraelGreville03} and \cite{CampbellMeyer09}). Considerations which follows
are described in terms of $\{1\}$-inverse of matrices.

\medskip
This section of the paper is organized as follows: In subsection {\em 2.1.} we
represent a new form of condition for the consistency of the matrix equation
(\ref{AXBC}). An extension of Penrose's theorem related to the general solution of the
matrix equation (\ref{AXBC}) is given in subsection {\em 2.2.} Namely, we represent
the formula of general solution of the matrix equation (\ref{AXBC}) if any particular
solution $X_{0}$ is known. In subsection {\em 2.3.} we give a form of particular
solution $X_{0}$ such that the formula of general solution of the matrix equation
(\ref{AXBC}), which is given in subsection {\em 2.2.}, is reproductive. The main
results of this paper are obtained in subsections {\em 2.1.}$-${\em 2.3.} and
additionally in subsection {\em 2.4.} we give two applications of the
concept~of~repro\-ductivity to some matrix systems which are in relation to the matrix
equation~(\ref{AXBC}).

\bigskip
\textbf{\textit{2.1.}}
Let $A\!\in\!\mathbb{C}_{a}^{m \times n}$, $B\!\in \!\mathbb{C}_{b}^{p \times q}$ and
$C\!\in\!\mathbb{C}^{m \times q}$. The matrix  $A \!\in \!\mathbb{C}_{a}^{m \times n}$
has $a$ linearly independent rows and $a$ linearly independent columns.
Let $T_{A_{r}}$ be a $m \times m$ permutation matrix such that multiplying the matrix
$A$ by the matrix $T_{A_{r}}$ on the left, we can permute the rows of the matrix $A$
and let $T_{A_{c}}$ be a $n \times n$ permutation matrix such that multiplying the
matrix $A$ by the matrix $T_{A_{c}}$ on the right, we can permute the columns of the
matrix $A$. Then, for the matrix $A$ there are permutation matrices $T_{A_{r}}$ and
$T_{A_{c}}$ such that the matrix
\begin{equation}
\label{AKAPA}
\widehat{A}=T_{A_{r}}AT_{A_{c}}
\end{equation}
has linearly independent rows and linearly independent columns \textit{at the first
$a$ positions}. Analogously, for the matrix $B$ there are permutation matrices
$T_{B_{r}}$ and $T_{B_{c}}$ such that
the matrix
\begin{equation}
\label{BKAPA}
\widehat{B}=T_{B_{r}}BT_{B_{c}}
\vspace*{-0.7 mm}
\end{equation}
has linearly independent rows and linearly independent columns \textit{at the first
$b$ positions.}

\smallskip
\noindent
The considerations which follow are valid for any choice of matrices $ T_{A_{r}}$,
$T_{A_{c}}$, $T_{B_{r}}$ and $T_{B_{c}}$ such that $\widehat{A}$ has linearly
independent rows and linearly independent columns at the first $a$ positions and
$\widehat{B}$ has linearly independent rows and linearly independent columns at the
first $b$ positions.~Let
\vspace*{-1.5 mm}
\begin{equation}
\label{CKAPA}
\widehat{C}=T_{A_{r}}CT_{B_{c}}.
\end{equation}
Next, let for the matrices $A$ and $B$ regular matrices $Q_{1}, \, P_{1}$ and $Q_{2},
\, P_{2}$ be determined such that the following equalities are true:
\vspace*{-1.0 mm}
\begin{equation}
\label{EAEB}
Q_{1}AP_{1}=
E_{a}=
\left[
\begin{array}{c|c}
I_{a} & 0\\\hline
0 & 0
\end {array}
\right]
\qquad \mbox{and} \qquad
Q_{2}BP_{2}=E_{b}=\left[
\begin{array}{c|c}
I_{b} & 0\\\hline
0 & 0
\end {array}
\right]
\end{equation}
\vspace*{-1.0 mm}
i.e.
\vspace*{-1.0 mm}
\begin{equation}
\label{AB}
A=Q_{1}^{-1}E_{a}P_{1}^{-1}
\qquad \mbox{and} \qquad
B=Q_{2}^{-1}E_{b}P_{2}^{-1}.
\end{equation}

\break

\noindent
Then, from (\ref{AKAPA}), (\ref{BKAPA}) and (\ref{AB}) we  get  that
$$
\qquad \widehat{A}=T_{A_{r}}Q_{1}^{-1}E_{a}P_{1}^{-1}T_{A_{c}}
\qquad \mbox{and} \qquad
\widehat{B}=T_{B_{r}}Q_{2}^{-1}E_{b}P_{2}^{-1}T_{B_{c}} \qquad
$$
i.e.
$$
\quad \widehat{A}=(Q_{1}T_{A_{r}}^{-1})^{-1}E_{a}(T_{A_{c}}^{-1}P_{1})^{-1}
\quad \mbox{and}  \quad
\widehat{B}=(Q_{2}T_{B_{r}}^{-1})^{-1}E_{b}(T_{B_{c}}^{-1}P_{2})^{-1}.\quad
$$
If we introduce the following notations:
\begin{equation}
\label{QPKAPA}
\widehat{Q_{1}} = Q_{1}T_{A_{r}}^{-1}, \;\; \widehat{P_{1}} = T_{A_{c}}^{-1}P_{1}
\quad \mbox{and} \quad
\widehat{Q_{2}} = Q_{2}T_{B_{r}}^{-1}, \;\; \widehat{P_{2}} = T_{B_{c}}^{-1}P_{2}
\end{equation}
we get that
\begin{equation}
\label{ABKAPA}
\widehat{A}=\widehat{Q_{1}}^{-1}E_{a} \widehat{P_{1}}^{-1} \qquad \mbox {and}  \qquad
\widehat{B}=\widehat{Q_{2}}^{-1}E_{b}\widehat{P_{2}}^{-1}.
\end{equation}
Considering Rohde's general form of $\{1\}$-inverses $A^{(1)}$ and $B^{(1)}$:
\begin{equation}
\label{ABJEDAN}
A^{(1)}=
P_{1}
\left[
\begin{array}{c|c}
I_{a} & X_{1}\\\hline
X_{2} & X_{3}
\end {array}
\right]
Q_{1}
\qquad \mbox{and} \qquad
B^{(1)}=
P_{2}
\left[
\begin{array}{c|c}
I_{b} & Y_{1}\\\hline
Y_{2} & Y_{3}
\end{array}
\right]
Q_{2} \, ,
\end{equation}
where  $X_{1}$, $X_{2}$, $X_{3}$ and $Y_{1}$, $Y_{2}$, $Y_{3}$ are arbitrary
matrices of suitable sizes, we obtain that:
$$
\qquad \widehat{A}^{(1)}=
\widehat{P_{1}}
\left[
\begin{array}{c|c}
I_{a} & X_{1}\\\hline
X_{2} & X_{3}
\end{array}
\right]
\widehat{Q_{1}}
\qquad \mbox{and} \qquad
\widehat{B}^{(1)}=
\widehat{P_{2}}
\left[
\begin{array}{c|c}
I_{b} & Y_{1}\\\hline
Y_{2} & Y_{3}
\end {array}
\right]\widehat{Q_{2}} \qquad
$$
i.e.
\vspace*{-0.5 mm}
\begin{equation}
\label{AAJEDANKAPA}
\widehat{A}\widehat{A}^{(1)}=
\widehat{Q_{1}}^{-1}E_{a}\widehat{P_{1}}^{-1}\widehat{P_{1}}
\left[
\begin{array}{c|c}
I_{a} & X_{1}\\\hline
X_{2} & X_{3}
\end{array}
\right]
\widehat{Q_{1}}=
\widehat{Q_{1}}^{-1}
\left[
\begin{array}{c|c}
I_{a} & X_{1}\\\hline
0 & 0
\end{array}
\right]
\widehat{Q_{1}}
\end{equation}

\smallskip
\noindent
and
\vspace*{-0.5 mm}
\begin{equation}
\label{BJEDANBKAPA}
\widehat{B}^{(1)}\widehat{B}=
\widehat{P_{2}}
\left[
\begin{array}{c|c}
I_{b} & Y_{1}\\\hline
Y_{2} & Y_{3}
\end {array}
\right]
\widehat{Q_{2}}\widehat{Q_{2}}^{-1}E_{b}\widehat{P_{2}}^{-1}=
\widehat{P_{2}}
\left[
\begin{array}{c|c}
I_{b} & 0\\\hline
Y_{2} & 0
\end {array}
\right]\widehat{P_{2}}^{-1}.
\end{equation}

\smallskip
\noindent
As we mentioned, the matrix $\widehat{A}$ has linearly independent rows and linearly
independent columns at the first $a$ positions and the matrix $\widehat{B}$ has
linearly independent rows and linearly independent columns at the first $b$ positions.

\smallskip
\noindent
Let
\vspace*{-0.5 mm}
\begin{equation}
\label{VRSTEAKAPA}
\widehat{A}_{i \rightarrow}= \sum_{l=1}^{a}\alpha_{i,l} \widehat{A}_{l \rightarrow},
\quad i=a+1,...,m \, ,
\end{equation}

\begin{equation}
\label{KOLONEAKAPA}
\widehat{A}_{\downarrow j}= \sum_{k=1}^{a}\alpha'_{k,j}\widehat{A}_{\downarrow k},
\quad j=a+1,...,n \, ,
\end{equation}

\noindent
and

\begin{equation}
\label{VRSTEBKAPA}
\widehat{B}_{i \rightarrow}=\sum_{l=1}^{b}\beta'_{i,l}\widehat{B}_{l \rightarrow},
\qquad i=b+1,...,p \, ,
\end{equation}

\begin{equation}
\label{KOLONEBKAPA}
\widehat{B}_{\downarrow j}=\sum_{k=1}^{b}\beta_{k,j}\widehat{B}_{\downarrow k},\qquad
j=b+1,...,q \, ;
\end{equation}

\break

\noindent
for some scalars $\alpha_{i,l}$, $\alpha'_{k,j}$ and $\beta'_{i,l}$, $\beta_{k,j}$. As
we know, the matrices $\widehat{Q_{1}}$, $\widehat{P_{1}}$ and $\widehat{Q_{2}}$,
$\widehat{P_{2}}$ are not uniquely determined, but we shall use, without loss of
generality, their following forms:
\begin{equation}
\label{OQPKAPAA}
\widehat{Q_{1}}=
\left[
\begin{array}{c|c}
I_{a} & 0 \\\hline
L_{1} & I_{m-a}
\end{array}
\right], \qquad
\widehat{P_{1}}=
\left[
\begin{array}{c|c}
W^{-1} & L'_{1} \\\hline
0 & I_{n-a}
\end{array}
\right]
\end{equation}
and
\begin{equation}
\label{OQPKAPAB}
\widehat{Q_{2}}=
\left[
\begin{array}{c|c}
U^{-1}& 0 \\\hline
L'_{2} & I_{p-b}
\end{array}
\right], \qquad
\widehat{P_{2}}=
\left[
\begin{array}{c|c}
I_{b} & L_{2} \\\hline
0 & I_{q-b}
\end{array}
\right]
\end{equation}

\medskip
\noindent
for

\medskip

\centerline{
$\quad L_{1}=\left[
\begin{array}{ccc}
-\alpha_{a+1,1}&...& -\alpha_{a+1,a}\\
 .&...& .\\
 .&...& .\\
 .&...& .\\
 -\alpha_{m,1}&...& -\alpha_{m,a}\\
\end{array}
\right]$,
\,
$L'_{1}=\left[
\begin{array}{ccc}
-\alpha'_{1,a+1}&...& -\alpha'_{1,n}\\
 .&...& .\\
 .&...& .\\
 .&...& .\\
 -\alpha'_{a,a+1}&...& -\alpha'_{a,n}\\
\end{array}
\right]$,}

\smallskip

\bigskip
\centerline{
$\quad L'_{2}=\left[
\begin{array}{ccc}
-\beta'_{b+1,1}&...& -\beta'_{b+1,b}\\
 .&...& .\\
 .&...& .\\
 .&...& .\\
 -\beta'_{p,1}&...& -\beta'_{p,b}\\
\end{array}
\right]$,
\,
$L_{2}=\left[
\begin{array}{ccc}
-\beta_{1,b+1}&...& -\beta_{1,q}\\
 .&...& .\\
 .&...& .\\
 .&...& .\\
 -\beta_{b,b+1}&...& -\beta_{b,q}\\
\end{array}
\right]$}

\medskip
\noindent
and where $W$ is a $a \times a $ submatrix of $\widehat{A}$ such that
$\widehat{A}=\!\left[
\begin{array}{c|c}
W & \widehat{A}_{2} \\ \hline
\widehat{A}_{3} & \widehat{A}_{4}
\end{array}
\right]$
and $U$ is a $b \times b$ submatrix of $\widehat{B}$ such that
$\widehat{B}=\!\left[
\begin{array}{c|c}
U & \widehat{B}_{2} \\ \hline
\widehat{B}_{3} & \widehat{B}_{4}
\end{array}
\right]$.

\bigskip
\noindent
Let us emphasize that the following statement is true.

\begin{lemma}
\label{L21}
Let $A\!\in\!\mathbb{C}_{a}^{m \times n}$, $B\!\in \!\mathbb{C}_{b}^{p \times q}$,
$C\!\in\!\mathbb{C}^{m \times q}$.
Suppose that $\widehat{A}$ and $\widehat{B}$ are determined by (\ref{AKAPA}) and
(\ref{BKAPA}). Then, the conditions
\vspace*{-0.5 mm}
\begin{equation}
\label{UKABC}
AA^{(1)}CB^{(1)}B=C
\end{equation}
\vspace*{-0.7 mm}
\noindent
and
\vspace*{-0.7 mm}
\begin{equation}
\label{UKABCKAPA}
\widehat{A}\widehat{A}^{(1)}\widehat{C}\widehat{B}^{(1)}\widehat{B}=\widehat{C}
\end{equation}
\vspace*{-0.5 mm}
\noindent
are equivalent.
\end{lemma}

\smallskip
\noindent
{\bf Proof.}$\,$The following equivalences are true
\mbox{$\!AA^{(1)}CB^{(1)}\!B \!=\! C$}
$\!\Longleftrightarrow\!$
\mbox{$T_{A_{r}}^{-1}\!\widehat{A}\!\;\!\widehat{A}^{(1)}\!\;\!T_{A_{r}}\!C\!\;\!T_{B_{c}}\!\widehat{B}^{(1)}\!\widehat{B}\!\;\!T_{B_{c}}^{-1} \!\!=\! C$}
$\!\Longleftrightarrow\!$
\mbox{$\widehat{A}\widehat{A}^{(1)}T_{A_{r}}\!CT_{B_{c}}\widehat{B}^{(1)}\widehat{B} \!=\! T_{A_{r}}\!CT_{B_{c}}$}
$\!\Longleftrightarrow\!$
\mbox{$\widehat{A}\widehat{A}^{(1)}\widehat{C}\widehat{B}^{(1)}\widehat{B} \!=\! \widehat{C}$}.
$\diamondsuit$

\bigskip
\noindent
Let us remark that (\ref{UKABC}) is {\em Penrose's condition of consistency}
for the matrix equation (\ref{AXBC}), \cite{Penrose55}, and if the matrix equation
(\ref{AXBC}) is consistent, then the formulas of general solution are given in Theorem
\ref{T22}.$\;$and$\;$\ref{T23}.$\;$In the following statement we give a condition which is
equivalent to Penrose's condition of consistency for the matrix equation
(\ref{AXBC}). So, we can use this new condition to test the consistency of the matrix
equation (\ref{AXBC}).

\medskip
\begin{theorem}
\label{T21}
Let $A\!\in \!\mathbb{C}_{a}^{m \times n}$,$B\!\in\!\mathbb{C}_{b}^{p \times q}$,
$C\!\in \!\mathbb{C}^{m \times q}$. Suppose that $\widehat{A}$ and $\widehat{B}$ are
determined by (\ref{AKAPA}) and (\ref{BKAPA}) and that
(\ref{VRSTEAKAPA})--(\ref{KOLONEBKAPA}) are satisfied. Then, the
condition (\ref{UKABC}) is true for any choice of  \textit{\{}$\!$1\textit{\}}-inverses
$A^{(1)}$ and $B^{(1)}$ iff
\begin{equation}
\label{MatrixC}
\mbox{
$\widehat{C}$=$\left[
\begin{array}{cccccc}
c_{1,1}\!&\!...\!&\!c_{1,b}\!&\!\mbox{\scriptsize
$\displaystyle\sum_{k=1}^{b}$}\beta_{k,b+1}c_{1,k}\!&\!...\!&\mbox{\scriptsize
$\displaystyle\sum_{k=1}^{b}$}\beta_{k,q}c_{1,k}\\[-0.75 ex]
 .\!&\!...\!&.&.&\!...\!&.\\[-0.75 ex]
 .\!&\!...\!&.&.&\!...\!&.\\[-0.75 ex]
 .\!&\!...\!&.&.&\!...\!&.\\[-0.75 ex]
 c_{a,1}\!&\!...\!&c_{a,b}&\mbox{\scriptsize
 $\displaystyle\sum_{k=1}^{b}$}\beta_{k,b+1}c_{a,k}&\!...\!&\mbox{\scriptsize
 $\displaystyle\sum_{k=1}^{b}$}\beta_{k,q}c_{a,k} \\[1.75 ex]
\mbox{\scriptsize $\displaystyle\sum_{l=1}^{a}$}\alpha_{a+1,l}c_{l,1}\!&\!...\!&
\mbox{\scriptsize $\displaystyle\sum_{l=1}^{a}$}\alpha_{a+1,l}c_{l,b}\!&
\begin{small}\mbox{\scriptsize $\displaystyle\sum_{l=1}^{a}$}\mbox{\scriptsize
$\displaystyle\sum_{k=1}^{b}$}\alpha_{a+1,l}\beta_{k,b+1}c_{l,k}\end{small} &\!...\!&
\begin{small}\mbox{\scriptsize $\displaystyle\sum_{l=1}^{a}$}\mbox{\scriptsize
$\displaystyle\sum_{k=1}^{b}$} \alpha_{a+1,l}\beta_{k,q}c_{l,k}\end{small}\\[-0.75 ex]
 .&\!...\!& .& .&\!...\!& .\\[-0.75 ex]
 .&\!...\!& .& .&\!...\!& .\\[-0.75 ex]
 .&\!...\!& .& .&\!...\!& .\\[-0.75 ex]
\mbox{\scriptsize $\displaystyle\sum_{l=1}^{a}$}\alpha_{m,l}c_{l,1}&\!...\!&
\mbox{\scriptsize $\displaystyle\sum_{k=1}^{b}$}\alpha_{m,l}c_{l,b}&
\begin{small}\mbox{\scriptsize $\displaystyle\sum_{l=1}^{a}$}\mbox{\scriptsize
$\displaystyle\sum_{k=1}^{b}$}\alpha_{m,l}\beta_{k,b+1}c_{l,k} \end{small}&\!...\!&
\begin{small} \mbox{\scriptsize $\displaystyle\sum_{l=1}^{a}$}\mbox{\scriptsize
$\displaystyle\sum_{k=1}^{b}$}\alpha_{m,l}\beta_{k,q}c_{l,k} \end{small}
\end {array}
\right]$,}
\end{equation}
where $c_{i,j}$ are arbitrary elements of $\mathbb{C}$.
\end{theorem}

\noindent
{\bf Proof.} {\boldmath $(\Longrightarrow)$}: Suppose that the condition (\ref{UKABC}) is valid for any choice of $\{1\}$-inverses
$A^{(1)}$ and $B^{(1)}.$ Based on Lemma \ref{L21}. the condition (\ref{UKABCKAPA}) is also valid. Then, considering the equalities
(\ref{AAJEDANKAPA}) and (\ref{BJEDANBKAPA}), we get the following  equality
\vspace*{-0.5 mm}
$$
\widehat{Q_{1}}^{-1}\left[ \begin{array}{c|c}
I_{a} & X_{1}\\\hline
0 & 0
\end {array}
\right]\widehat{Q_{1}}\widehat{C}\widehat{P_{2}}
\left[
\begin{array}{c|c}
I_{b} & 0\\\hline
Y_{2} & 0
\end {array} \right] \widehat{P_{2}}^{-1}
=
\widehat{C}.
$$
\noindent
By multiplying the previous equality by $ \widehat{Q_{1}}$ on the left and by $\widehat{P_{2}}$
on the right we get
\begin{equation}
\label{QCPKAPA}
\left[
\begin{array}{c|c}
I_{a} & X_{1}\\\hline
0 & 0
\end {array}
\right]
\widehat{Q_{1}}\widehat{C}\widehat{P_{2}}
\left[
\begin{array}{c|c}
I_{b} & 0\\\hline
Y_{2} & 0
\end {array}
\right]=
\widehat{Q_{1}}\widehat{C}\widehat{P_{2}}.
\end{equation}

\smallskip
\noindent
Suppose that

\medskip
\noindent
$$
\widehat{C}=
\left[\begin{array}{ccc}
c_{1,1}&...& c_{1,q}\\
 .&...& .\\
 .&...& .\\
 .&...& .\\
 c_{m,1}&...& c_{m,q}\\
\end {array}\right] \! .
$$

\medskip
\noindent
We are going to show that $\widehat{C}$ has the form (\ref{MatrixC}). Let
\noindent
\begin{equation}
\label{EF}
E=\widehat{Q_{1}}\widehat{C}\widehat{P_{2}} \qquad \mbox {and} \qquad
F=\left[ \begin{array} {c|c}
I_{a} & X_{1}\\\hline
0 & 0
\end {array}
\right]E\left[
\begin{array}{c|c}
I_{b} & 0 \\\hline
Y_{2} & 0
\end {array}
\right].
\end{equation}
From (\ref{OQPKAPAA}) and (\ref{OQPKAPAB}) we obtain that

\medskip
\noindent
for  $i=1,...,a, j=1,...,b$  \hspace*{19 mm} $E_{i,j}=c_{i,j},$

\vspace*{-1 mm}
\noindent
for $i=1,...,a, j=b+1,...,q$  \hspace*{12 mm} $E_{i,j}=c_{i,j}-\displaystyle\sum_{k=1}^{b}\beta_{k,j}c_{i,k},$

\vspace*{-1 mm}
\noindent
for $i=a+1,...,m, j=1,...,b$ \hspace*{10 mm}
$E_{i,j}=c_{i,j}-\displaystyle\sum_{l=1}^{a}\alpha_{i,l}c_{l,j},$

\noindent
for $i=a+1,...,m, j=b+1,...,q$\hspace*{5 mm}
$E_{i,j}=c_{i,j}-\displaystyle\sum_{l=1}^{a}\alpha_{i,l}c_{l,j}
        -\displaystyle\sum_{k=1}^{b}\beta_{k,j}(c_{i,k}-\displaystyle\sum_{l=1}^{a}\alpha_{i,l}c_{l,k})$

\noindent
and

\smallskip
\noindent
for $ i=1,...,a, j=1,...,b$\hspace*{21 mm} $F_{i,j}=c_{i,j}+
\displaystyle \sum_{\overline{l}=a+1}^{m}x_{i,
\overline{l}}(c_{\overline{l},j}-\displaystyle\sum_{l=1}^{a}\alpha_{\overline{l},l}c_{l,j})$

\hspace*{60 mm} $+\displaystyle\sum_{\overline{k}=b+1}^{q}y_{\overline{k},j}[ c_{i,
\overline{k}}-\displaystyle\sum_{k=1}^{b}\beta_{k,\overline{k}}c_{i,k}$

\hspace*{60 mm} $+\displaystyle\sum_{\overline{l}=a+1}^{m}x_{i,\overline{l}}
\lbrace c_{\overline{l},\overline{k}}-\displaystyle\sum_{l=1}^{a}\alpha_{\overline{l},l}c_{l,\overline{k}}$

\hspace*{60 mm} $-\displaystyle\sum_{k=1}^{b}\beta_{k,\overline{k}}
(c_{\overline{l},k}-\displaystyle\sum_{l=1}^{a}\alpha_{\overline{l},l}c_{l,k})\rbrace ] ,$

\smallskip
\noindent
for $i=1,...,a, j=b+1,...,q$  \hspace*{14 mm} $F_{i,j}=0,$

\smallskip
\noindent
for $ i=a+1,...,m, j=1,...,b$ \hspace*{12.5 mm}  $F_{i,j}=0,$

\smallskip
\noindent
for $ i=a+1,...,m, j=b+1,...,q$ \hspace*{5.75 mm}  $F_{i,j}=0.$

\medskip
\noindent
Finally, from (\ref{QCPKAPA}) and (\ref{EF}) i.e. $E=F$  we get that

\medskip
\noindent
for $ i=1,...,a, j=1,...,b$ \hspace*{18.5 mm} $c_{i,j}$ are arbitrary elements
of $\mathbb{C},$

\vspace*{-0.5 mm}
\noindent
for $ i=a+1,...,m, j=1,...,b$ \hspace*{10.5 mm}
$c_{i,j}=\displaystyle\sum_{l=1}^{a}\alpha_{i,l}c_{l,j},$

\vspace*{-0.5 mm}
\noindent
for $ i=1,...,a, j=b+1,...,q$ \hspace*{11.5 mm}
$c_{i,j}=\displaystyle\sum_{k=1}^{b}\beta_{k,j}c_{i,k},$

\vspace*{-0.5 mm}
\noindent
for $ i=a+1,...,m, j=b+1,...,q$ \hspace*{4 mm}
$c_{i,j}=\displaystyle\sum_{l=1}^{a}\sum_{k=1}^{b}\alpha_{i,l}\beta_{k,j}c_{l,k}.$

\smallskip
\noindent
{\boldmath $(\Longleftarrow)$}: Suppose that the matrix $\widehat{C}$ has the form (\ref{MatrixC}). Then,
\begin{equation}
\label{OQCPKAPA}
\qquad \widehat{Q_{1}}\widehat{C} \widehat{P_{2}}=...=\left[
\begin{array}{cccccc}
c_{1,1}&...&c_{1,b}&0&...& 0   \\[-0.5 ex]
 .&...& .&.&...& .             \\[-0.5 ex]
 .&...& .&.&...& .             \\[-0.5 ex]
 .&...& .&.&...& .             \\[-0.5 ex]
 c_{a,1}&...&c_{a,b}&0&...& 0  \\[-0.0 ex]
0&...& 0&0 &...& 0             \\[-0.5 ex]
 .&...& .& .&...& .            \\[-0.5 ex]
 .&...& .& .&...& .            \\[-0.5 ex]
 .&...& .& .&...& .            \\[-0.5 ex]
0&...& 0&0 &...& 0
\end{array}
\right]
\end{equation}

\break

\noindent
and
\begin{equation}
\label{POQCPKAPA}
\left[ \begin{array}{c|c}
I_{a} & X_{1}\\\hline
0 & 0
\end {array}
\right]\widehat{Q_{1}}\widehat{C} \widehat{P_{2}}\left[
\begin{array}{c|c}
I_{b}& 0 \\\hline
Y_{2}& 0
\end {array}
\right]=...=\left[
\begin{array}{cccccc}
c_{1,1}&...&c_{1,b}&0&...& 0\\
 .&...& .&.&...& .\\
 .&...& .&.&...& .\\
 .&...& .&.&...& .\\
 c_{a,1}&...&c_{a,b}&0&...& 0\\
0&...& 0&0 &...& 0\\
 .&...& .& .&...& .\\
 .&...& .& .&...& .\\
 .&...& .& .&...& .\\
0&...& 0&0 &...& 0
\end{array}
\right].
\end{equation}
From (\ref{OQCPKAPA}) and (\ref{POQCPKAPA}) we conclude that

$$
\left[ \begin{array}{c|c}
I_{a} & X_{1}\\\hline
0 & 0
\end {array}
\right]\widehat{Q_{1}}\widehat{C}
\widehat{P_{2}}\left[
\begin{array}{c|c}
I_{b}& 0 \\\hline
Y_{2}& 0
\end {array}
\right]= \widehat{Q_{1}}
\widehat{C}
\widehat{P_{2}}.
$$

\noindent
By multiplying the previous equality by $\widehat{Q_{1}}^{-1} $ on the left and by
$\widehat{P_{2}}^{-1}$ on the right we obtain the following equality:

\vspace*{-1.5 mm}

\begin{equation}
\label{NOCKAPA}
\widehat{Q_{1}}^{-1}
\left[
\begin{array}{c|c}
I_{a} & X_{1}\\\hline
0 & 0
\end{array}
\right]
\widehat{Q_{1}}\widehat{C}\widehat{P_{2}}
\left[
\begin{array}{c|c}
I_{b}& 0 \\\hline
Y_{2}& 0
\end {array}
\right]\widehat{P_{2}}^{-1}=\widehat{C}.
\end{equation}

\smallskip
\noindent
From (\ref{NOCKAPA}), considering the equalities (\ref{AAJEDANKAPA}) and
(\ref{BJEDANBKAPA}), we see that the condition (\ref{UKABCKAPA}) is true. Based on
Lemma \ref{L21}. we conclude that the condition~(\ref{UKABC})~is~true.
$\diamondsuit$

\begin{remark}
Let us remark that the general form of matrix $C$, such that the matrix equation
(\ref{AXBC}) is consistent, always exists. The matrix equation (\ref{AXBC}) is
consistent for an arbitrary matrix $C$ iff a matrix $A$ has full row rank and a matrix
$B$ has full column rank (see also Exercises 10.50 from \cite{KMAbadirJRMagnus}).
\end{remark}

\begin{remark}
In the paper \cite{Wang04} author considered some forms which are equivalent to
Penrose's condition of consistency for matrix equation (\ref{AXBC}).
\end{remark}

\bigskip
\noindent
The application of Theorem \ref{T21} will be illustrated by the following examples.

\begin{example}
\label{E21}
Let be given the following matrices:

\bigskip
\begin{center}
$A= \left[
\begin{array}{rrr}
 0 &  0 &  0\\
 1 & -3 &  2\\
 2 &  1 & -1\\
-1 & -4 &  3\\
 3 & -2 &  1
\end{array}
\right]$
and
$\;B= \left[
\begin{array}{rrrrr}
0 & 1 & 2 & 3 \!&\!-1\\
0 & 3 & 1 & 4 \!&\! 2\\
0 & 4 & 1 & 5 \!&\! 3\\
0 & 2 & 3 & 5 \!&\!-1
\end {array}
\right].$
\end{center}

\medskip
\noindent
Then, $rank(A)$=2, $rank(B)$=2 and for

\medskip
\begin{center}
$T_{A_{r}}=
\left[
\begin{array}{ccccc}
0&1&0&0&0\\
0&0&1&0&0\\
0&0&0&1&0\\
0&0&0&0&1\\
1&0&0&0&0
\end{array}
\right]$
and
$\;
T_{B_{c}}= \left[
\begin{array}{ccccc}
0&0&1&0&0\\
1&0&0&0&0\\
0&1&0&0&0\\
0&0&0&1&0\\
0&0&0&0&1
\end{array}
\right]$
\end{center}
we get that
$$
\widehat{A}=T_{A_{r}}A
=...=
\left[
\begin{array}{rrr}
 1 & -3 &  2\\
 2 &  1 & -1\\
-1 & -4 &  3\\
 3 & -2 &  1\\
 0 &  0 &  0
\end {array}
\right]\!,
\;\;
\widehat{B}=BT_{B_{c}}
=... =
\left[
\begin{array}{rrrrr}
 1 & 2 & 0 & 3 & -1\\
 3 & 1 & 0 & 4 &  2\\
 4 & 1 & 0 & 5 &  3\\
 2 & 3 & 0 & 5 & -1
\end{array}
\right]\!.
$$

\smallskip
\noindent
Therefore,

\smallskip
\begin{center}
$ \widehat{A}_{3\rightarrow}=\widehat{A}_{1\rightarrow}-\widehat{A}_{2\rightarrow}$,
\hspace*{1 mm}
$ \widehat{A}_{4\rightarrow}=\widehat{A}_{1\rightarrow}+\widehat{A}_{2\rightarrow}$,
\hspace*{1 mm}
$\widehat{A}_{5\rightarrow}=0\widehat{A}_{1\rightarrow}+0\widehat{A}_{2\rightarrow}$
\end{center}
\, and \,
\begin{center}
$\widehat{B}_{\downarrow 3}=0\widehat{B}_{\downarrow 1}+0\widehat{B}_{\downarrow 2}$,
\hspace*{1 mm}
$ \widehat{B}_{\downarrow 4}=\widehat{B}_{\downarrow 1}+\widehat{B}_{\downarrow 2}$,
\hspace*{1 mm}
$\widehat{B}_{\downarrow 5}=\widehat{B}_{\downarrow 1}-\widehat{B}_{\downarrow 2}.$
\end{center}
From this we get that
\begin{center}
$ \alpha_{3,1}=1,$ $\alpha_{3,2}=-1,$  $ \alpha_{4,1}=1,$ $ \alpha_{4,2}=1,$
$\alpha_{5,1}=0,$  $\alpha_{5,2}=0 $
\end{center}
and
\begin{center}
$ \beta_{1,3}=0, $ $ \beta_{2,3}=0 ,$  $ \beta_{1,4}=1, $  $ \beta_{2,4}=1, $ $ \beta_{1,5}=1, $
$ \beta_{2,5}=-1. $
\end{center}
Based on Theorem \ref{T21}. each matrix $\widehat{C}$ which has the following form
$$
\widehat{C}
=
\mbox{
\scriptsize $\left[
\begin{array}{ccccc}
c_{1,1}&c_{1,2}&0&c_{1,1}+c_{1,2}& c_{1,1}-c_{1,2}\\
c_{2,1}&c_{2,2}&0&c_{2,1}+c_{2,2}&c_{2,1}-c_{2,2}\\
c_{1,1}-c_{2,1}&c_{1,2}-c_{2,2}&0&c_{1,1}-c_{2,1}+c_{1,2}-c_{2,2}&c_{1,1}-c_{2,1}-
c_{1,2}+c_{2,2}\\
c_{1,1}+c_{2,1}&c_{1,2}+c_{2,2}&0&c_{1,1}+c_{2,1}+c_{1,2}+c_{2,2}&c_{1,1}+c_{2,1}-
c_{1,2}-c_{2,2}\\
0&0& 0& 0& 0
\end{array}
\right]$
}
$$
satisfies the condition (\ref{UKABCKAPA}). From that we conclude that each matrix
$C=T_{A_{r}}^{-1}\widehat{C}T_{B_{c}}^{-1}$ which has the following form
$$
C
=
\mbox{
\scriptsize $\left[
\begin{array}{ccccc}
0&0&0&0& 0\\
0&c_{1,1}&c_{1,2}&c_{1,1}+c_{1,2}& c_{1,1}-c_{1,2}\\
0&c_{2,1}&c_{2,2}&c_{2,1}+c_{2,2}& c_{2,1}-c_{2,2}\\
0&c_{1,1}-c_{2,1}&c_{1,2}-c_{2,2}&c_{1,1}-c_{2,1}+c_{1,2}-c_{2,2}&c_{1,1}-c_{2,1}-
c_{1,2}+c_{2,2}\\
0&c_{1,1}+c_{2,1 }&c_{1,2}+c_{2,2}&c_{1,1}+c_{2,1 }+c_{1,2}+c_{2,2}&c_{1,1}+c_{2,1}-
c_{1,2}-c_{2,2}
\end {array}
\right]$
}
$$
satisfies the condition (\ref{UKABC}).
$\blacklozenge$
\end{example}

\begin{example}
\label{E22}
Let $A$ and $B$ be the matrices as in Example \ref{E21}. and

\vspace*{7 mm}

a) $\;
C=\left[
\begin{array}{rrrrr}
0 & 0 & 0 & 0 & 0 \\
0 & 1 & 0 & 1 & 1 \\
0 &-2 & 2 & 0 &-4 \\
0 & 3 &-2 & 1 & 5 \\
0 &-1 & 2 & 1 &-3
\end {array}
\right]\!,$
\hspace*{2 mm}
b) $\;
C=\left[
\begin{array}{rrrrr}
0 & 0 & 0 & 0 & 0 \\
0 & 1 & 0 & 1 & 1 \\
0 &-2 & 2 & 0 & 4 \\
0 & 3 &-2 & 1 & 5 \\
0 &-1 & 2 & 1 &-3
\end {array}
\right]\!.$

\vspace*{7 mm}

\noindent
If we compare the matrix $C$ from a) and from b) with the general form of matrix
$C$ which satisfies the condition (\ref{UKABC}) (see Example \ref{E21}.) we see that
the matrix $C$ from a) satisfies the condition (\ref{UKABC}) and the matrix $C$
from~b) does not satisfy the condition (\ref{UKABC}). Therefore, the matrix
equation~(\ref{AXBC})~is~consistent for the matrix $C$ from a), but it is not
consistent for the matrix $C$ from~b).~$\blacklozenge $
\end{example}

\bigskip
\textbf{\textit{2.2.}}
Recall that the matrix equation $AXB=C$ is marked with (\ref{AXBC}) for $A\!\in\!
\mathbb{C}_{a}^{m \times n}$, \mbox{$B\!\in \!\mathbb{C}_{b}^{p \times q}$}, $C\!\in\!
\mathbb{C}^{m \times q}$. Methods for solving the consistent matrix equation
(\ref{AXBC}) are considered in the book \cite{Cullis1913} (Chapter~X). In the paper
\cite{Penrose55} R. Penrose proved the following theorem related to the matrix
equation (\ref{AXBC}).

\begin{theorem}
\label{T22}
The matrix equation (\ref{AXBC}) is consistent iff for some choice of \textit{\{}$\!
$1\textit{\}}-inverses $A^{(1)}$ and
$B^{(1)}$ of the  matrices $A$ and $B$ the condition (\ref{UKABC}) is true. The
general solution of the matrix equation
(\ref{AXBC}) is given by the formula

\smallskip

\begin{equation}
\label{ORR}
X=f(Y)=A^{(1)}CB^{(1)}+Y-A^{(1)}AYBB^{(1)},
\end{equation}

\smallskip

\noindent
where $Y\!\in\!\mathbb{C}^{n \times p}$ is an arbitrary matrix.
\end{theorem}

\begin{remark}
\label{R21}
If the matrix equation (\ref{AXBC}) is consistent, the equivalence

\smallskip

\begin{equation}
AXB=C
\;\Longleftrightarrow\;
X=f(X)=X-A^{(1)}(AXB-C)B^{(1)}
\end{equation}

\smallskip

\noindent
is true. Therefore, the starting equation is equivalent to some reproductive equation.
Based~on~Theorem \ref{T12}. we can also conclude that (\ref{ORR}) is the general
solution of the matrix~equation~(\ref{AXBC}).
\end{remark}

\smallskip
\noindent
In this paper we give a simple extension of Theorem \ref{T22}.

\begin{theorem}
\label{T23}
If $X_{0}$ is any particular solution of the matrix equation (\ref{AXBC}), the
general solution of the matrix equation (\ref{AXBC}) is given by the formula
\begin{equation}
\label{OR}
X=g(Y)=X_{0}+Y-A^{(1)}AYBB^{(1)},
\end{equation}
where $Y\!\in\!\mathbb{C}^{n \times p}$ is an arbitrary matrix. The function $g$
satisfies the condition of reproductivity (\ref{UR}) iff $X_{0}=A^{(1)}CB^{(1)}$.
\end{theorem}

\medskip
\noindent
{\bf Proof.}
It is easily to see that the solution of the matrix equation (\ref{AXBC}) is given by
(\ref{OR}). On the contrary, let $ X $ is any solution of the matrix equation
(\ref{AXBC}), then
\begin{eqnarray*}
X \!&\!=\!&\! X-A^{(1)}CB^{(1)}+A^{(1)}CB^{(1)}          \\
  \!&\!=\!&\! X-A^{(1)}AXBB^{(1)}+A^{(1)}AX_{0}BB^{(1)}  \\
  \!&\!=\!&\! X-A^{(1)}A(X-X_{0})B B^{(1)}               \\
  \!&\!=\!&\! X_{0}+(X-X_{0})-A^{(1)}A(X-X_{0})BB^{(1)}  \\
  \!&\!=\!&\! X_{0}+Y-A^{(1)}AYBB^{(1)}=g(Y) \, ,
\end{eqnarray*}
where $Y=X-X_{0}$. From this we see that every solution $X$ of the matrix equation
(\ref{AXBC}) can be represented in the form (\ref{OR}). Based on the following matrix
equality:
$$
g^{2}(Y)
=
g(Y)+(X_{0}-A^{(1)}CB^{(1)})
$$
we see that the function $g$ satisfies the condition (\ref{UR}) iff
$X_{0}=A^{(1)}CB^{(1)}$.
$\diamondsuit$
\begin{remark}
Using the previous theorem and the appropriate choice of particular solution $X_{0}$
we can obtain the general solutions for different cases of the matrix equation
(\ref{AXBC}). It was considered in the papers \cite{Haveric83} and \cite{Presic63}.
\end{remark}

\bigskip
\noindent
The general solution (\ref{OR}) of the matrix equation (\ref{AXBC}) is reproductive
iff $X_{0}=A^{(1)}CB^{(1)}$. Therefore, {\em Penrose's general solution} (\ref{ORR})
of the matrix equation (\ref{AXBC}) is the reproductive solution. If the condition
(\ref{UKABC}) is not true, the matrix equation (\ref{AXBC}) is solved
approximately as described in the paper \cite{Penrose55} and books
\cite{KMAbadirJRMagnus}, \cite{Ben-IsraelGreville03} and \cite{CampbellMeyer09}.

\bigskip
\textbf{\textit{2.3.}}
Using the obtained form of matrix $\widehat{C}$ we obtain the form of particular
solution $X_{0}$ of the matrix equation (\ref{AXBC}) such that the general solution
(\ref{OR}) of the matrix equation (\ref{AXBC}) is reproductive.

\begin{theorem}
\label{T24}
Let $X_{0}$ any particular solution of the matrix equation (\ref{AXBC}). The
general solution (\ref{OR}) of  the  matrix  equation (\ref{AXBC}) is reproductive iff
\begin{equation}
\label{Form_X_0}
X_{0}=P_{1}
\left[
\begin{array}{c|c}
C_{1} & C_{1}Y_{1}\\\hline
X_{2}C_{1} & X_{2}C_{1}Y_{1}
\end {array}
\right]Q_{2}
\end{equation}
where $P_{1}$, $Q_{2}$, $ X_{2} $, $ Y_{1}$ are the  matrices  from (\ref{ABJEDAN})
and $C_{1} $ is the submatrix of the matrix  $\widehat{C}$ and it has the following
form:
$$
C_{1}=\left[
\begin{array}{ccc}
c_{1,1}&...& c_{1,b}\\
 .&...& .\\
 .&...& .\\
 .&...& .\\
 c_{a,1}&...& c_{a,b}\\
\end {array}
\right],
$$
where $c_{i,j}$ are some elements of $\mathbb{C}$.
\end{theorem}

\medskip
\noindent
{\bf Proof.}
For general $\{1\}$-inverses $A^{(1)}$ and $B^{(1)}$ the statement follows from
Theorem \ref{T23}. because
$$
\begin{array}{rcl}
X_{0}
&\!\!\!=\!\!\!&
A^{(1)}C B^{(1)} \mathop{=} \limits_{(\ref{ABJEDAN})} P_{1}
\left[
\begin{array}{c|c}
I_{a} & X_{1}\\\hline
X_{2} & X_{3}
\end {array}
\right]
Q_{1}CP_{2}
\left[
\begin{array}{c|c}
I_{b} & Y_{1}\\\hline
Y_{2} & Y_{3}
\end {array}
\right]Q_{2}                                                                       \\[1.5 ex]
&\!\!\!\mathop{=} \limits_{(\ref{CKAPA})}\!\!\!&
P_{1}
\left[
\begin{array}{c|c}
I_{a} & X_{1}\\\hline
X_{2} & X_{3}
\end {array}
\right]Q_{1}T_{A}^{-1}\widehat{C}T_{B}^{-1} P_{2}
\left[
\begin{array}{c|c}
I_{b} & Y_{1}\\\hline
Y_{2} & Y_{3}
\end {array}
\right]Q_{2}                                                                       \\[1.5 ex]
&\!\!\!\mathop{=} \limits_{(\ref{QPKAPA})}\!\!\!&
P_{1}
\left[
\begin{array}{c|c}
I_{a} & X_{1}\\\hline
X_{2} & X_{3}
\end {array}
\right]\widehat{Q_{1}}\widehat{C}\widehat{ P_{2}}
\left[
\begin{array}{c|c}
I_{b} & Y_{1}\\\hline
Y_{2} & Y_{3}
\end {array}
\right]Q_{2}                                                                       \\[1.5 ex]
&\!\!\!\mathop{=} \limits_{(\ref{OQCPKAPA})}\!\!\!&
P_{1}
\left[
\begin{array}{c|c}
I_{a} & X_{1}\\\hline
X_{2} & X_{3}
\end {array}
\right]
\left[\begin{array}{c|c}
C_{1} & 0\\\hline
0 & 0
\end {array}
\right]
\left[
\begin{array}{c|c}
I_{b} & Y_{1}\\\hline
Y_{2} & Y_{3}
\end {array}
\right]Q_{2}                                                                       \\[1.5 ex]
&\!\!\!=\!\!\!&
P_{1}
\left[
\begin{array}{c|c}
C_{1} & C_{1}Y_{1}\\\hline
X_{2}C_{1} & X_{2}C_{1}Y_{1}
\end {array}
\right]Q_{2} \,.\quad \diamondsuit
\end{array}
$$
\begin{corollary}
\label{X_0_Full_rank}
\textit{\textbf{(i)}} If a matrix $A$ has full row rank and a matrix $B$ has full column rank, then parameters
from submatrices $X_2$ and $Y_1$ don't exist and don't appear in the matrix $X_{0}$ of~form~(\ref{Form_X_0}).
\textit{\textbf{(ii)}} If either a matrix $A$ has full row rank or a matrix $B$ has full column rank, then
the matrix $X_0$ of form (\ref{Form_X_0}) has the structure of an affine linear space with parameters from
either submatrix~$X_2$~or~submatrix~$Y_1$, respectively. \textit{\textbf{(iii)}} If a matrix~$A$ doesn't have
full row rank and a matrix~$B$ doesn't have full column rank, then the matrix $X_0$ of form (\ref{Form_X_0})
doesn't have the structure of an affine linear space relative to parameters from submatrices $X_2$ and $Y_1$.
\end{corollary}
\begin{remark}
\label{R23}
In the paper \cite{MalesevicRadicic11} authors proved that there is a matrix equation
(\ref{AXBC}) and its particular solution $X_{1}$ such that
$ X_{1}\neq A^{(1)}CB^{(1)}$ for any choice of \textit{\{}$\!$1\textit{\}}-inverses
$A^{(1)}$ and $B^{(1)}$.
\end{remark}
\begin{remark}
\label{R24}
According to Theorem VI, pp.$\,$345-346, from \cite{Cullis1913}, it is possible to extract $a \!\cdot\! b$ parameters in a matrix $Y$ such that,
these parameters are expressed, in the solution (\ref{OR}), as a non-homogeneous linear functions of the other $n\!\cdot\! p - a \!\cdot\! b$
independent parameters.
\end{remark}

\bigskip
\textbf{\textit{2.4.}}
In this part of the paper we analysed two applications the concept of reproductivity
on some matrix systems which are in relation to the matrix equation (\ref{AXBC}).
\begin{application}
In \cite{Penrose55} R.~Penrose studied a matrix system
\begin{equation}
\label{S1}
(\ref{S1}a) \quad AX=B
\qquad \wedge \qquad
(\ref{S1}b) \quad XD=E,
\qquad
\end{equation}
where $A,$ $B,$ $D$ and $E$ are given complex matrices corresponding dimensions.
He proved that
\begin{equation}
X_{1}=A^{(1)}B+ED^{(1)}-A^{(1)}AED^{(1)}
\end{equation}
is one common solution of the matrix equations (\ref{S1}a) and (\ref{S1}b) if $AE=BD$ and the matrix equations
(\ref{S1}a) and (\ref{S1}b) are consistent.

\bigskip
\noindent
In \cite{Ben-IsraelGreville03} A. Ben-Israel and  T.N.E. Greville proved that the matrix equations (\ref{S1}a) and
(\ref{S1}b) have a common solution iff each equation separately has a solution and  $AE=BD$. Also, they proved
that if $X_{0}$ is any common solution of the matrix equations (\ref{S1}a) and (\ref{S1}b), the general solution  of
the matrix  system (\ref{S1}) is given by the  formula
\begin{equation}
X=g(Y)=X_{0}+(I-A^{(1)}A)Y(I-DD^{(1)}),
\end{equation}
where $ Y $ is an  arbitrary  matrix  corresponding dimensions.

\break

\noindent
We will prove that if the matrix system (\ref{S1}) is consistent, the general reproductive solution is
given by the  formula
\begin{equation}
\label{ORRS1}
X=f(Y)=A^{(1)}B+ED^{(1)}-A^{(1)}AED^{(1)}+(I-A^{(1)}A)Y(I-DD^{(1)}),
\end{equation}
where $Y$ is an arbitrary matrix corresponding dimensions.

\medskip
\noindent
If the matrix system (\ref{S1}) is consistent, the following equivalence is true
\begin{equation}
\label{App_1}
{\big (}AX \!=\!B \; \; \wedge \; \; XD \!=\! E{\big )}
\;\Longleftrightarrow\;
X=f(X).
\end{equation}
The direct implication of (\ref{App_1}) follows by implications (see Remark \ref{R21}. in~the~subsection~2.1):

\vspace*{-2.0 mm}

$$
AX \!=\! B \;\Longrightarrow\; X \!=\! f_{1}(X) \!=\! A^{(1)}B \!+\! X \!-\! A^{(1)}AX \, ,
$$
$$
XD \!=\! E \;\Longrightarrow\; X \!=\! f_{2}(X) \!=\! ED^{(1)} \!+\! X \!-\!XDD^{(1)} \, ,
$$
$$
AXD \!=\! BD \!=\! AE \;\Longrightarrow\; X \!=\! f_{3}(X) \!=\! A^{(1)}AED^{(1)} \!+\! X \!-\! A^{(1)}AXDD^{(1)} \, .
$$

\smallskip
\noindent
From the previous implications we can conclude
$$
{\big (} AX \!=\! B \; \; \wedge \; \; XD \!=\! E {\big )}
\;\Longrightarrow\;
X \!=\! f(X) \!=\! f_{1}(X) \!+\! f_{2}(X) \!-\! f_{3}(X).
$$

\noindent
The reverse implication of (\ref{App_1}) is trivial. Notice that the function $f$ is reproductive. Therefore,
if the matrix system (\ref{S1}) is consistent, it is equivalent to the reproductive matrix equation
\mbox{$X\!=\!f(X)$}. Based on Theorem \ref{T12}. we conclude that $X\!=\!f(Y)$ is the general
reproductive solution of the~matrix system (\ref{S1}). If there is a particular solution $X_0$ of
the matrix system (\ref{S1}) so that $X_0 \!\neq\! X_1$, then $X \!=\! g(Y)$ is the general non-reproductive
solution. At the end of these application let us remark that equality $X \!=\! g(Y\!-\!X_0) \!=\! f(Y)$ also
represents one simple proof of the Statement 1 from \cite{Haveric84}.
\end{application}

\begin{application}
\label{A22}Let $A \in \mathbb{C}^{n \times n}$ be a singular matrix.
In this section we consider a matrix system
\begin{equation}
\label{Sys2}
AXA=A \quad \wedge \quad AX=XA.
\end{equation}

\medskip
\noindent
The consistency of the matrix system (\ref{Sys2}) is determined by Theorem 1 in \cite{Keckic85} (see also \cite{Ruski_Gruv_84_85}~and~\cite{Keckic97}).
Let $\bar{A}$ is commutative \textit{\{}1\textit{\}}-inverse, \cite{Keckic85}. Based on the reproductivity, we give a new proof that the formula from
\cite{Keckic85}:
\begin{equation}
\label{fSolA2}
X
=
f(Y)
=
\bar{A} A \bar{A} + Y - \bar{A} A Y - Y A \bar{A} + \bar{A} A Y A \bar{A},
\end{equation}
where $Y$ is an arbitrary matrix corresponding dimensions, represents the general solution of the consistent matrix system (\ref{Sys2}).

\medskip
\noindent
Namely, if the matrix system (\ref{Sys2}) is consistent, the equivalence
\begin{equation}
\label{App_2}
{\big (} AXA=A \; \; \wedge \; \; AX=XA {\big )}
\;\Longleftrightarrow\;
X=f(X)
\end{equation}
is true. The direct implication of (\ref{App_2}) is based on the following simple matrix equalities:
$$
\bar{A} \!\!\mathop{\underbrace{\!\!A\!\!}}\limits_{\mbox{\scriptsize $(=\!A\!X\!A)$}}\!\! \bar{A}
=
\bar{A}AX \!\!\mathop{\underbrace{A\bar{A}}} \limits_{\mbox{\scriptsize $(=\bar{A}\!A)$}}
=
\bar{A} \!\!\mathop{\underbrace{AX}}\limits_{\mbox{\scriptsize $(=XA)$}}\!\! \bar{A}A
=
\bar{A}X \!\!\mathop{\underbrace{A\bar{A}A}}\limits_{\mbox{\scriptsize $(=A)$}}
=
\bar{A} \!\!\mathop{\underbrace{XA}}\limits_{\mbox{\scriptsize $(=AX)$}}
=
\bar{A}AX
$$
and
$$
\mathop{\underbrace{\bar{A} A}}\limits_{\mbox{\scriptsize $(=A\bar{A})$}}\! X A \bar{A}
=
A\bar{A} \!\mathop{\underbrace{X A}}\limits_{\mbox{\scriptsize $(=AX)$}}\! \bar{A}
=
\mathop{\underbrace{A\bar{A} A}}\limits_{\mbox{\scriptsize $(=A)$}}\! X\bar{A}
=
\mathop{\underbrace{A X}}\limits_{\mbox{\scriptsize $(=XA)$}}\! \bar{A}
=
XA\bar{A}.
$$

\break

\noindent
From this we get that $X \!=\! X \!+\! \bar{A} A \bar{A} \!-\! \bar{A}AX \!+ \! \bar{A} A X A \bar{A} \!-\! X A \bar{A} \!=\! f(X)$.
The reverse implication of (\ref{App_2}) is trivial. Notice that the function $f$ is reproductive. Therefore, if the matrix system
(\ref{Sys2}) is consistent, it is equivalent to the reproductive matrix equation $X \!=\! f(X)$. Based on Theorem \ref{T12}. we conclude
that~\mbox{$X \!=\! f(Y)$} is the general reproductive solution of the matrix system (\ref{Sys2}). If $X_{0}$ is any solution of the matrix
system (\ref{Sys2}), the formula
\begin{equation}
\label{gSolA2}
X = g(Y) = X_{0} + Y - \bar{A} A Y - Y A \bar{A} + \bar{A} A Y A \bar{A},
\end{equation}
also determines a form of the general solution of the matrix system (\ref{Sys2})
because the equality \mbox{$g(Y) \!=\! f(X_0 + Y)$} is true.
If there is a particular solution $X_0$ of the matrix system (\ref{Sys2}) such that
$X_0 \!\neq\! \bar{A} A \bar{A}$, then $X \!=\!g(Y)$ is the
general non-reproductive solution. Additional applications of the concept of reproductivity for some matrix equations and systems were considered in the paper
\cite{MalesevicRadicic12}.
\end{application}

\medskip
\noindent
{\bf Acknowledgment.} Research is partially supported by the Ministry of Science and Education of the Republic of Serbia,
Grant No. ON 174032.


\bigskip

\begin{thebibliography}{20}

\bibitem{KMAbadirJRMagnus} {\sc K.M. Abadir and J.R. Magnus}, {\em Matrix Algebra}, Econometric exercises,
Volume \textbf{1}, Cambridge, 2005.

\bibitem{Bankovic79} {\sc D. Bankovi\' c}, {\em On general and reproductive solutions of arbitrary equations},
Publications  de l'insti\-tut math\' e\-mati\-que, Nouvelle serie, tome \textbf{26} (40), Beograd 1979, 31$\;$-$\,$33.

\bibitem{Bankovic95} {\sc D. Bankovi\' c}, {\em All solutions of finite equations}, Discrete Mathematics
Vol. {\bf 137} (1-3), 1995, 1$\;$-$\,$6.

\bibitem{Bankovic97} {\sc D. Bankovi\' c}, {\em General reproductive solutions of Postian equations}, Discrete Mathematics
Vol. {\bf 169} (1-3), 1997, 163$\;$-$\,$168.

\bibitem{Bankovic02} {\sc D. Bankovi\' c}, {\em All general solutions of Pre\v si\' c's equation},
Facta universitatis, Ser. Math. Inform. Vol. {\bf 17}, Ni\v s 2002, 1$\;$-$\,$4.

\bibitem{Bankovic11} {\sc D. Bankovi\' c}, {\em General Solutions of System of Finite Equations},
Scientific Publications of the State University of Novi Pazar Ser. A: Appl. Math. Inform. and
Mech. vol. \textbf{3}, 2 (2011), 117$\;$-$\,$121.

\bibitem{Ben-IsraelGreville03} {\sc A. Ben-Israel and T.N.E. Greville}, {\em Generalized Inverses: Theory
and Applications}, Springer, 2003.

\bibitem{Bozic75} {\sc  M.$\;$Bo\v zi\' c}, {\em A Note On Reproductive Solutions}, Publications
de l'insti\-tut math\' e\-mati\-que, Nouvelle serie, tome \textbf{19} (33), Beograd 1975, 33$\;$-$\,$35.
({\tt http:/$\!$/publications.mi.sanu.ac.rs/})

\bibitem{CampbellMeyer09} {\sc S.L. Campbell and C.D. Meyer}, {\em Generalized Inverses of Linear
Transformations}, Society for Industrial and Applied Mathematics, 2009.

\bibitem{Cullis1913} {\sc C.E. Cullis}, {\em Matrices and determinoids$\,$-$\!\,$Volume $I$}, Cambridge, University Press~Publ.~1913.
({\tt http:/$\!$/archive.org/details/matricesdetermin01cull})

\break

\bibitem{Cvetkovic06} {\sc D.S. Cvetkovi\' c-Ili\' c}, {\em The reflexive solutions of the matrix equation
$AXB=C$}, Comp. Math. Appl., \textbf{51} (2006), 897$\;$-$\,$902.

\bibitem{Cvetkovic07}
{\sc D.S. Cvetkovi\' c-Ili\' c, A. Daji\' c, and J.J. Koliha}, {\em Positive and real-positive solutions to the equation $axa^{\ast}\!=\!c$
in $C^{\ast}-\,$algebra}, Linear \& Multilinear algebra, \textbf{55}, (6) (2007), 535$\;$-$\,$543.

\bibitem{Cvetkovic08} {\sc D.S. Cvetkovi\' c-Ili\' c}, {\em Re-nnd solutions of the matrix equation
$AXB=C$}, Journal of the Australian Mathematical Society, \textbf{84} (2008), 63$\;$-$\,$72.

\bibitem{ADajicJJKoliha} {\sc A. Daji\' c and J.J. Koliha}, {\em Equations $ax=c$ and $xb=d$ in rings and rings with involution
with applications to Hilbert space operators}, Lin. Alg. and its Appl. {\bf 429} (2008) 1779$\;$-$\,$1809.

\bibitem{Ruski_Gruv_84_85} {\sc L.D. Dobryakov}, {\em Commuting generalized inverse matrices}, Mathematical
Notes, Volume \textbf{36}, Number 1, 500$\;$-$\,$504, 1985 (Translated from Matematicheskie Zametki, Vol.
\textbf{36}, No. 1, 17$\,$-$\,$23, 1984.).

\bibitem{V. Harizanov} {\sc V. Harizanov}, {\em On the functional equation $f \phi f \!=\! f$},
Publications de l'institut  math\' e\-matique, Nouvelle serie, tome \textbf{29} (43), Beograd 1981, 61$\,$-$\,$64.

\bibitem{Haveric83} {\sc M. Haveri\' c}, {\em Formulae for general reproductive solutions of certain matrix
equations}, Publications de l'institut  math\' ematique, Nouvelle serie, tome \textbf{34} (48), Beograd 1983, 81$\,$-$\,$84.

\bibitem{Haveric84} {\sc M. Haveri\' c}, {\em On solutions of a matrix equations system $AX=B$ and $XD=E$},
Matemati\v cki Vesnik \textbf{36} (1), Beograd 1984, 11$\,$-$\,$16.

\bibitem{Keckic82} {\sc J.D. Ke\v cki\' c}, {\em Reproductivity of some equations of analysis I},
Publications de l'insti\-tut math\' e\-mati\-que, Nouvelle serie, tome \textbf{31}(45), Beograd 1982, 73$\;$-$\,$81.

\bibitem{Keckic83} {\sc J.D. Ke\v cki\' c}, {\em Reproductivity of some equations of analysis II},
Publications de l'insti\-tut math\' e\-mati\-que, Nouvelle serie, tome \textbf{33}(47), Beograd 1983, 109$\;$-118.

\bibitem{Keckic85} {\sc J.D. Ke\v cki\' c}, {\em Commutative weak generalized inverses of a square matrix
and some related matrix equations}, Publications  de l'insti\-tut math\' e\-mati\-que, Nouvelle
serie, tome \textbf{38} (52), Beograd 1985, 39$\;$-$\,$44.

\bibitem{Keckic86} {\sc J.D. Ke\v cki\' c}, {\em On some generalized inverses of matrices and some linear
matrix  eguations}, Publications  de l'insti\-tut math\' e\-mati\-que, Nouvelle serie, tome \textbf{45} (59),
Beograd 1989, 57$\,$-$\,$63.

\bibitem{Keckic97} {\sc J.D. Ke\v cki\' c}, {\em Some remarks on possible generalized inverses in semigroups},
Publications  de l'insti\-tut math\' e\-mati\-que, Nouvelle serie, tome \textbf{61} (75), Beograd 1997, 33$\;$-$\,$40.

\bibitem{KeckicPresic97} {\sc J.D. Ke\v cki\' c and S.B. Pre\v si\' c}, {\em Reproductivity - A general approach to equations},
Facta universitatis, Ser. Math. Inform. Vol. {\bf 12}, Ni\v s 1997, 157$\;$-$\,$184.

\bibitem{MalesevicRadicic11} {\sc B. Male\v sevi\' c and B. Radi\v ci\' c}, {\em Non-reproductive and
reproductive solutions~of~some matrix equations}, Proceedings of the International conference {\em Mathematical and
Informational Technologies, MIT$\,$-$\,$2011}, Vrnja\v cka Banja, Serbia, 2011, 246$\;$-$\,$251.
({\tt http:/$\!$/mit.rs/})

\bibitem{MalesevicRadicic12} {\sc B. Male\v sevi\' c and B. Radi\v ci\' c}, {\em Some considerations of matrix equations
using the concept of reproductivity}, Kragujevac Journal of Mathematics, \textbf{36}(1) (2012), 151$\;$-$\,$161.

\break

\bibitem{Penrose55} {\sc R. Penrose}, {\em A generalized inverses for matrices}, Math. Proc. Cambridge
Philos. Soc. \textbf{51}(1955), 406$\;$-$\,$413.

\bibitem{PresicC.R.Acad.Sci.Paris63} {\sc S.B. Pre\v si\' c}, {\em Methode de resolution d'une classe d'equations fonctionnelles lineaires},
Comptes rendus de l'Acad\' emie des Sciences Paris, {\bf 257} (1963), 2224$\;$-$\,$2226.

\bibitem{Presic63} {\sc S.B. Pre\v si\' c}, {\em Certaines  \' equations matricielles}, Publ. Elektrotehn.
Fak. Ser. Mat.-Fiz., $\mathit{N}^{\underline{o}}$~121, Beograd 1963. ({\tt http:/$\!$/pefmath.etf.rs/})

\bibitem{Presic68} {\sc S.B. Pre\v si\' c}, {\em Une classe d'\' equations  matricielles et l'\' equation
fonctionnelle $f^{2}\!=\!f$}, Publications de l'institut  math\' ematique, Nouvelle serie, tome \textbf{8}
(22), Beograd 1968, 143$\;$-$\,$148.

\bibitem{PresicC.R.Acad.Sci.Paris71} {\sc S.B. Pre\v si\' c}, {\em  Une methode de resolution des equations dont toutes les solutions
appartiennent a un ensemble fini donne}, Comptes rendus de l'Acad\' emie des Sciences Paris, {\bf 272} (1971), \mbox{654$\;$-$\,$657}.

\bibitem{Presic72} {\sc S.B. Pre\v si\' c}, {\em Ein Satz \" Uber Reproduktive L\" osungen},
Publications de l'in\-sti\-tut math\' e\-ma\-ti\-que, Nouvelle serie, tome \textbf{14} (28),
Beograd 1972, 133$\;$-136.

\bibitem{Presic88} {\sc S.B. Pre\v si\' c}, {\em All reproductive solutions of finite equations},
Publications de l'in\-sti\-tut math\' e\-ma\-ti\-que, Nouvelle serie, tome \textbf{44} (58),
Beograd 1988, 3$\,$-7.

\bibitem{Presic00} {\sc S.B. Pre\v si\' c}, {\em A generalization of the notion of reproductivity},
Publications de l'insti\-tut math\' e\-mati\-que, Nouvelle serie, tome \textbf{67} (81), Beograd 2000, 76$\;$-$\,$84.

\bibitem{Rohde64} {\sc C.A. Rohde}, {\em  Contribution to the theory, computation and application of
generalized inverses}, Doctoral  dissertation, University of North Carolina at Releigh, May 1964.

\bibitem{Rudeanu78} {\sc S. Rudeanu}, {\em On general solutions of arbitrary equations},
Publications  de l'insti\-tut math\' e\-mati\-que, Nouvelle serie, tome \textbf{24} (38), Beograd 1978, 143$\;$-$\,$145.

\bibitem{Rudeanu98} {\sc S. Rudeanu}, {\em On general and reproductive solutions of finite equations},
Publications  de l'insti\-tut math\' e\-mati\-que, Nouvelle serie, tome \textbf{63} (77), Beograd 1998, 26$\;$-$\,$30.

\bibitem{Rudeanu01} {\sc S. Rudeanu}, {\em Lattice Functions and Equations}, Springer, 2001.

\bibitem{Tian10} {\sc Y. Tian}, {\em On additive decompositions of solutions of the matrix equation $AXB=C$},
Calcolo, Vol. \textbf{47} (4), 2010, 193$\;$-$\,$209.

\bibitem{Tian12} {\sc Y. Tian}, {\em On Additive Decomposition of the Hermitian Solution of the Matrix Equation $AXA^{\ast}=B$},
Mediterranean Journal of Mathematics, \textbf{9} (2012), 47$\;$-$\,$60.

\bibitem{Wang04} {\sc Q$\,$-W. Wang}, {\em A system of matrix equations and a linear matrix equation over arbitrary regular rings with identity},
Lin. Alg. and its Appl. \textbf{384}, (2004), 43$\;$-$\,$54.


\bibitem{Tribute01}
{\sc A. Krape\v z (editor)}, {\it A tribute to S.$\,$B. Pre\v si\' c$\,:$} {\sl Papers Celebrating his 65$\,$-$\,$th~\mbox{\sl Birthday}},
Mathematical Institute of the Serbian Academy of Sciences and Arts, Belgrade, publ. 2001.
({\tt http:/$\!$/elibrary.matf.bg.ac.rs/handle/123456789/448})

\break

\end{thebibliography}
\end{document}